\documentclass[reqno]{amsart}

\usepackage{xypic}
\input xy
\xyoption{all}
\usepackage{epsfig}
\usepackage{color}
\usepackage{amsthm}
\usepackage{amssymb}
\usepackage{amsmath}
\usepackage{amscd}
\usepackage{amsopn}

\newcommand{\labell}[1] {\label{#1}}

\numberwithin{equation}{section}
\newtheorem {Theorem}{Theorem} 
\numberwithin{Theorem}{section}

\newtheorem {Lemma}[Theorem]    {Lemma}         
\newtheorem {Proposition}[Theorem]{Proposition}  
\theoremstyle{definition}

\theoremstyle{remark}
\newtheorem{Remark}[Theorem]{Remark}
\newtheorem{Example}[Theorem]{Example}

\newtheorem {Corollary}[Theorem]{Corollary}  

\expandafter\chardef\csname pre amssym.def at\endcsname=\the\catcode`\@
\catcode`\@=11
\def\undefine#1{\let#1\undefined}
\def\newsymbol#1#2#3#4#5{\let\next@\relax
 \ifnum#2=\@ne\let\next@\msafam@\else
 \ifnum#2=\tw@\let\next@\msbfam@\fi\fi
 \mathchardef#1="#3\next@#4#5}
\def\mathhexbox@#1#2#3{\relax
 \ifmmode\mathpalette{}{\m@th\mathchar"#1#2#3}%
 \else\leavevmode\hbox{$\m@th\mathchar"#1#2#3$}\fi}
\def\hexnumber@#1{\ifcase#1 0\or 1\or 2\or 3\or 4\or 5\or 6\or 7\or 8\or
 9\or A\or B\or C\or D\or E\or F\fi}

\font\teneufm=eufm10
\font\seveneufm=eufm7
\font\fiveeufm=eufm5
\newfam\eufmfam
\textfont\eufmfam=\teneufm
\scriptfont\eufmfam=\seveneufm
\scriptscriptfont\eufmfam=\fiveeufm

\catcode`\@=\csname pre amssym.def at\endcsname

\def    \eps    {\epsilon}

\newcommand{\CH}{{\mathcal H}}
\newcommand{\CA}{{\mathcal A}}
\newcommand{\CalD}{{\mathcal D}}

\newcommand{\CM}{{\mathcal M}}

\newcommand{\CS}{{\mathcal S}}

\newcommand{\supp}{{\mathit supp}\,}

\newcommand{\id}{{\mathit id}}
\newcommand{\const}{{\mathit const}}

\def    \C      {{\mathbb C}}
\def    \R      {{\mathbb R}}
\def    \reals      {{\mathbb R}}
\def    \Z      {{\mathbb Z}}

\def    \T      {{\mathbb T}}
\def    \CP     {{\mathbb C}{\mathbb P}}

\def    \12    {{\frac{1}{2}}}

\def    \area  {\operatorname{area}}

\def    \SH     {\operatorname{SH}}
\def    \LD     {\underline{\operatorname{D}}}
\def    \HF     {\operatorname{HF}}

\def    \ssminus        {\smallsetminus}

\def    \CW  {\operatorname{c_{\scriptscriptstyle{HZ}}^o}}
\def    \wsh {\operatorname{c_{hom}}}
\def    \CHZ  {\operatorname{c_{\scriptscriptstyle{HZ}}}}
\def    \CGR  {\operatorname{c_{\scriptscriptstyle{Gr}}}}
\def    \cf  {\operatorname{c}}
\def    \H  {\operatorname{\scriptscriptstyle{H}}}
\def    \CHHZ  {{\mathcal H}_{\operatorname{\scriptscriptstyle{HZ}}}}

\begin{document}


\setlength{\smallskipamount}{6pt}
\setlength{\medskipamount}{10pt}
\setlength{\bigskipamount}{16pt}





\title[The Weinstein conjecture]{The Weinstein conjecture 
and the theorems of nearby and almost existence}

\author[Viktor Ginzburg]{Viktor L. Ginzburg}

\dedicatory{\bigskip\large{Dedicated to Alan Weinstein on the occasion of 
his sixtieth birthday}}

\address{Department of Mathematics, UC Santa Cruz, 
Santa Cruz, CA 95064, USA}
\email{ginzburg@math.ucsc.edu}

\subjclass[2000]{53D40, 37J45}
\date{\today}
\thanks{The work is partially supported by the NSF and by the faculty
research funds of the University of California, Santa Cruz.}

\bigskip

\begin{abstract}

The Weinstein conjecture, as the general existence problem for
periodic orbits of Hamiltonian or Reeb flows, has been among 
the central questions in symplectic topology for over two decades
and its investigation has led
to understanding of some fundamental properties of Hamiltonian flows.
 
In this paper we survey some recently developed and
well-known methods of proving various particular cases of this
conjecture and the closely related almost existence theorem. 
We also examine differentiability and continuity properties of the
Hofer--Zehnder capacity function and relate these properties to the
features of the underlying Hamiltonian dynamics, e.g., to the period
growth.

\end{abstract}

\maketitle

\section{Introduction} 
Without exaggeration, one can say that 
Arnold's conjecture and the Weinstein conjecture have been the two  
problems determining the development of symplectic topology over the past
twenty years. The Weinstein conjecture, \cite{We:conj}, the problem we
are interested in here, concerns the existence of closed characteristics on
a compact hypersurface of contact type. To be more precise, consider
a regular compact contact type level of a Hamiltonian on a symplectic manifold.
Then, the Weinstein conjecture as we understand it now asserts that the
level must carry at least one periodic orbit of the Hamiltonian flow.

This conjecture was motivated by numerous results establishing the
existence of periodic orbits under various, often rather restrictive,
conditions on the level (e.g., convex or bounding a star-shaped domain);
see, e.g., \cite{Mo:orbits,Ra,We:lagr,We:orbits,We:convex}. 
The feature distinguishing the Weinstein conjecture from these 
results is that the requirement that the level has contact type is invariant
under symplectomorphisms (as is the assertion) while the hypotheses of the
earlier theorems are not. In its original form
for Hamiltonians on a linear symplectic space,
the Weinstein conjecture was proved by Viterbo, \cite{Vi}. Since then,
the conjecture has been established for many other symplectic manifolds
(sometimes under additional requirements on the level); see, e.g.,
\cite{FHV, HV1, HV2, LT, Lu, Lu:uniruled,Lu:toric,Vi:functors}, 
to mention just a few results.
Among these manifolds are products of complex projective 
spaces, manifolds of the form $P\times \C^n$ (where $P$ is compact symplectic 
and the product is given a split symplectic structure), and
sub-critical Stein manifolds.
However, in general, the conjecture is still neither proved nor disproved.
For example, the Weinstein conjecture in its full generality is open for 
cotangent bundles and $\T^4$. 

Starting from the Weinstein conjecture, one can go in two different directions.
One direction is to dispose of the Hamiltonian and the ambient symplectic
manifold and focus exclusively on the level of contact type. This naturally
leads to the question of whether a Reeb flow on a compact contact manifold
necessarily has a periodic orbit; see, e.g., \cite{Ho:icm}.  
Along these lines, Hofer, \cite{Ho:sphere}, proved the existence of periodic
orbits for the Reeb flow of a contact form on $S^3$ or 
on a closed three-dimensional contact manifold $M$ with
$\pi_2(M)\neq 0$ and also for the Reeb flow of an overtwisted contact form.
This 
approach interprets the conjecture as a question about the dynamics of Reeb
flows and leads to the notions of contact homology and
symplectic field theory, \cite{EGH,BEHWZ}. 

Another direction is to view the conjecture as a question about the dynamics
of Hamiltonian flows on the ambient symplectic manifold. This is
the perspective with which we are concerned here. More specifically, we focus
on such 
problems as whether and how the assumption that the level has contact type
can be relaxed, whether the existence of periodic orbits is typical, etc.
The contact type requirement cannot be dropped entirely: there exists a 
proper function on $\R^{2n}$ ($C^\infty$-smooth if $2n\geq 6$ and $C^2$-smooth
if $2n=4$) with a regular level carrying no periodic orbits; see
\cite{gi:seifert95,gi:seifert97,gg1,gg2,He:fax,He,ke2}. 
(Constructions of such functions are known as counterexamples to the
Hamiltonian Seifert conjecture.) Nevertheless, the existence
of periodic orbits on the level sets of a fixed Hamiltonian
is a generic phenomenon: almost all (in the sense
of measure theory) regular levels of a $C^2$-smooth proper function on
$\R^{2n}$ carry periodic orbits of the Hamiltonian flow. This result,
due to Hofer and Zehnder and to Struwe, is known as the almost existence
theorem; \cite{HZ:book,St}. The almost existence theorem 
holds for many, but not all, symplectic manifolds (sometimes under some
additional restrictions on the Hamiltonian, cf. \cite{gg3}). 
For example, the theorem holds for $\CP^n$, products 
$P\times \C^n$ with split symplectic structure (where $P$ is geometrically 
bounded), \cite{HV1,Lu,MDS,FS3}, sub-critical Stein manifolds, 
\cite{Ci,Ke:new,Lu}, and small neighborhoods of certain non-Lagrangian
submanifolds, \cite{Ke:new,Sc,FS3}. 
However, as has been observed by Zehnder, \cite{Ze}, the almost existence 
theorem fails for certain Hamiltonians on $\T^{2n}$ equipped
with an irrational symplectic 
structure (Zehnder's torus). In fact, for this Hamiltonian system, there is 
an interval of energy values without periodic orbits.
Moreover, this phenomenon is stable under
$C^k$-small ($k>2n$) perturbations of the Hamiltonian, \cite{He2,He3}. 
A sibling of the almost existence theorem is the  theorem of dense or 
nearby existence
which guarantees the existence of periodic orbits for a dense set of
energy levels or, equivalently, near a fixed level; 
see, e.g., \cite{CGK,FH,FS1,HZ:cap,HZ:book}. 
Both of these theorems imply the Weinstein conjecture and below
we will discuss the  relation between these theorems. 

Another reason to consider the almost existence and nearby existence theorems 
is that a broad class of Hamiltonian systems for which energy
levels fail to have contact type  naturally arises in classical mechanics. 
In this class are, for example, the systems describing the motion of a charge 
in a magnetic field, which we will discuss shortly; see also \cite{gi:survey}.
For such systems one can expect the almost existence theorem to hold and
periodic orbits to exist for all low energy values. However, periodic orbits
need not exist on all energy levels even when the symplectic form
is exact near a level (the horocycle flow, see Example \ref{ex:horocycle}). 
Furthermore, as has been pointed out by Kerman, \cite{ke},
the analysis of such systems is closely related to a
generalization of the Weinstein--Moser theorem, \cite{Mo:orbits,We:orbits}.

The Weinstein--Moser theorem asserts that a smooth function $H$ on $\R^{2n}$
attaining a non-degenerate minimum at the origin must have at least
$n$ distinct periodic orbits on every level near the origin,
\cite{Mo:orbits,We:orbits}. Let us now replace $\R^{2n}$ by an arbitrary
symplectic manifold $W$ and the non-degenerate minimum at one point
by a Morse--Bott non-degenerate minimum along a closed symplectic submanifold
$M\subset W$. Then, conjecturally, every level of $H$ near $M$ carries
at least one periodic orbit or even a number of periodic orbits. 
We will refer to this conjecture as the generalized Weinstein--Moser 
conjecture.
As a particular example motivating our interest in this question, consider
the standard kinetic energy $H$ on the cotangent bundle 
$W=T^*M\stackrel{\pi}{\to}M$, equipped with 
a twisted symplectic structure $\omega_0+\pi^*\Omega$, where $\omega_0$ is
the standard symplectic structure $dp\wedge dq$ and $\Omega$ is a closed
two-form on $M$. This system describes the motion of a charge on $M$ in the
magnetic field $\Omega$; see \cite{gi:survey}. When $\Omega$ is 
non-degenerate $M$ turns into a symplectic submanifold of $W$. 

The generalized Weinstein--Moser conjecture has been proved in a number
of particular cases, \cite{gk1,ke}, but in general the question is still
open. Recently, however, some progress has been made along the lines of the 
almost existence theorem. Namely, it has been shown that almost all
levels close to $M$ of a function $H$ attaining a minimum along $M$ 
carry contractible periodic orbits, \cite{CGK,gg3,mac2} -- this is
the relative (with respect to $M$) almost existence theorem for small
energy values. (Note that unless
$M$ is a Morse--Bott non-degenerate minimum of $H$, we cannot expect
such periodic orbits to exist on all levels near $M$, \cite{gg3}.) This
implies the almost existence theorem for low energy periodic orbits of
a charge in a non-degenerate magnetic field. Moreover, under suitable
additional hypothesis, the genuine, non-relative, almost existence
theorem holds near $M$, \cite{Ke:new}. In the setting
of magnetic fields, these results can be further refined: periodic 
orbits
must exist whenever $\Omega\neq 0$, \cite{FS3}. (We will briefly discuss
the proof of this result in Section \ref{sec:Hofer-appl}.) We also refer the
reader to \cite{gk2, Ma,pol2} for related results.

The nearby existence theorem is weaker and often easier to prove than
the almost existence theorem. The pattern has been that, in many
cases,
the nearby existence theorem was proved first and then followed by the
almost existence theorem; cf., e.g., \cite{HZ:viterbo} and \cite{HZ:book,St},
\cite{CGK} and \cite{gg3}. As we have pointed
out above, both theorems imply the Weinstein conjecture. Conversely, 
essentially every proof of a particular case of the Weinstein conjecture in
the Hamiltonian setting translates into a proof of either the nearby
existence or almost existence theorem, although it is not always easy to
establish which of these theorems is proved. As of today, almost existence
is verified in virtually all the cases where the nearby existence has been 
proved. Probable exceptions are some of the results from \cite{Vi:functors} 
on periodic orbits in cotangent bundles and, perhaps, the results of 
\cite{LT} and hence of \cite{Lu:uniruled,Lu:GW,Lu:toric}.

In the present paper we focus on the aspects of the Weinstein conjecture
related to global Hamiltonian dynamics, and hence our treatment of
the conjecture in the large is by no means comprehensive. For example, we do 
not even touch upon the Weinstein conjecture for contact manifolds. (The reader
interested in this conjecture is referred to the surveys 
\cite{Ho:icm,Ho:parc-city,Ho:gafa}
in addition to the references given above.) Furthermore, we do not
mention the fruitful connection between the Weinstein conjecture and 
Gromov--Witten invariants, although we do discuss the holomorphic curve
approach to the proof of the conjecture.

In Sections \ref{sec:nearby} and \ref{sec:almost} we outline methods of
proving the nearby existence and almost existence theorems and, 
hence, the Weinstein conjecture. These sections can be viewed as a brief
introduction to certain concepts of symplectic topology (action selectors,
constructions of symplectic capacities, Hofer's metric, symplectic homology, 
etc.), albeit
strictly focused on a specific task and not even mentioning  many aspects 
of the subject. However, these sections should not be taken as an 
introduction to symplectic topology in general. For example, we assume 
the reader's familiarity with Floer 
homology. Section \ref{sec:cap-fun} concerns some simple, but apparently
not present in the literature, properties of the Hofer--Zehnder capacity 
function. Although the paper is complemented by an extensive 
bibliography, the list of references contains only the papers immediately
relevant to our discussion. Inevitably, this list is incomplete and 
omits many important contributions to the subject, and our exposition
emphasizes the publications that have most influenced the author's thinking.

\subsection*{Conventions.}
In this paper, all manifolds are assumed to be without boundary. A 
symplectic manifold
$W$ will be called convex if $W$ is either closed (i.e., compact) or open and
convex at infinity. Here $W$ is said to be convex at infinity if there
exists: a hypersurface $\Sigma\subset W$ which separates $W$ into
one set with compact closure and another, $U$, with non-compact
closure; and a flow $\varphi_t$ (for $t\geq 0$) of symplectic dilations 
on $U$, which is transversal to $\partial U=\Sigma$.
Recall also that $(W,\omega)$ is symplectically aspherical if 
$c_1|_{\pi_2(W)}=\omega|_{\pi_2(W)}=0$. (In some instances, this condition
can be replaced by $\omega|_{\pi_2(W)}=0$.) We refer the reader to, e.g.,
\cite{AL,CGK} for a discussion of geometrically bounded symplectic
manifolds and
to \cite{HZ:book,McDS:book,McDS:bookj,polt:book} for a general introduction 
to symplectic topology.

Let us now fix the sign conventions in the definition of the action
functional. Let $H\in C^\infty(S^1\times W)$ and let
$x\colon S^1\to W$ be a contractible loop in $W$. Here we will use
the action functional $A_H$ defined as
\begin{equation*}
A_H(x) = - \int_{D^2} \bar{x}^* \omega +
\int_{S^1}H(t,x)\,dt,
\end{equation*}
where $\bar{x}\colon D^2\to W$ is a map of a disk, bounded by $x$. 

By the least action principle, 
contractible one-periodic orbits of $H$ are precisely the critical 
points of $A_H$. The action spectrum $\CS(H)$ of $H$ is the set of critical 
values of $A_H$, i.e., the collection of action values of $A_H$ on the 
contractible
periodic orbits. (Since  we are assuming that  $\omega|_{\pi_2(W)}=0$,
the action $A_H$ is single valued. Otherwise, one
has to pass to a suitable covering of the loop space and deal with
numerous other technical difficulties.) We set $H_t=H(t,\cdot)$, where
$t\in S^1$, and denote by $X_H$ the Hamiltonian vector field of $H$. 

\subsection*{Acknowledgments.} The author is deeply grateful to 
Paul Biran, Urs Frauenfelder, Ba\c sak G\"urel, Ely Kerman, Debra Lewis, 
Leonid Polterovich, Felix Schlenk, and the referee 
for their useful remarks and suggestions.

\section{The nearby existence theorems}
\labell{sec:nearby}
Nearby existence theorems are usually proved by Floer homological methods.
Below we discuss some of these methods, focusing on the key
ideas and omitting many (often quite non-trivial) technical details.

Throughout this section, all ambient symplectic manifolds 
$(W,\omega)$ are assumed to be convex and symplectically aspherical,
unless explicitly stated otherwise as is, for example, in Section
\ref{sec:symhom}.

The argument is particularly transparent for open manifolds which are
convex at infinity. Hence we consider this case first.

\subsection{The action selector method: convex at infinity open manifolds} 
\labell{sec:as-open}
In this section we focus on convex at infinity open manifolds and outline 
the approach to the proof of the nearby existence theorem utilizing the notion 
of an action selector. In 
our proof of the existence of an action selector and its properties, 
we mainly follow \cite{FS1,Sc}. However, many elements of our argument
are already contained in 
\cite{HZ:cap,HZ:book},  where the action selector is defined differently;
see also, e.g., \cite{Oh2,Vi:functions}.

\subsubsection{Action selectors and the nearby existence theorem}
Let $W$ be a convex at infinity open symplectic manifold. 
We denote by $\CH_U$, where $U$ is an open subset 
of $W$, the space of compactly supported smooth Hamiltonians 
$H\colon S^1\times U\to \R$.

For our purposes, it is convenient to adopt the following definition.
An action selector is a $C^0$-continuous function $\sigma\colon \CH_W\to \R$ 
such that\footnote{The properties (AS1)--(AS3) should not be taken
as an attempt to axiomatize the notion of an action selector (cf. \cite{FGSS})
--- such
axioms would almost certainly be different from (AS1)--(AS3) and would include
a version of the inequality \eqref{eq:subad} below.
Our objective is to list the properties of the action selector that
are most essential for the proof of the nearby existence theorem.}

\begin{enumerate}

\item[(AS1)] $\sigma(H)\in \CS(H)$,

\item[(AS2)] $\sigma$ is monotone, i.e., $H_1\leq H_2$ implies that 
$\sigma(H_1)\leq \sigma(H_2)$,  and

\item[(AS3)] $\sigma(H)=\max H$, when $H$ is a $C^2$-small function
on $W$ with a unique maximum.
\end{enumerate}

Let us point out some consequences of (AS1)--(AS3) which are essential
for what follows. First recall that the action spectrum $\CS(H)$ is
compact and nowhere dense; see, e.g., \cite{HZ:book,Sc}. 
(This elementary, but not
entirely trivial, fact can be thought of as a version of Sard's theorem
for $A_H$, cf. \cite{Zv}.) Furthermore, as is easy to see, 
$\CS(H)$ depends only on the time-one flow $\varphi_H$ of $H$ viewed
as an element of the universal cover $\CalD$ of the group of Hamiltonian 
symplectomorphisms; the homotopy to identity is given by $H$.
(In fact, when $W$ is symplectically aspherical, $\CS(H)$ depends only
on the time-one flow of $H$ regarded as an element of the group of
symplectomorphisms.)  Then
it follows from (AS1) and the continuity of $\sigma$ that $\sigma(H)$
also depends only on $\varphi_H\in \CalD$. Furthermore, 
\begin{equation}
\labell{eq:positivity}
\sigma(H)>0\text{ when $H$ is non-negative and non-constant,}
\end{equation}
by (AS2) and (AS3) and, by (AS1) and (AS2), $\sigma(H)=0$ when $H\leq 0$.

The nearby existence theorem for functions on $U\subset W$ can be 
(and often is) easily derived from the existence of an action selector $\sigma$
which meets one additional requirement that $\sigma$ is \emph{a priori} bounded
on $\CH_U$. To state the result, it is convenient to introduce the
notion of a shell or thickening of a hypersurface.
A shell $\{\Sigma_s\}$ in $U$ is an embedding
$\Sigma\times (-\eps,\eps)\hookrightarrow U$, where $\Sigma$ is a closed
manifold such that $\dim \Sigma =\dim U-1$. We identify
$\Sigma\times (-\eps,\eps)$ with its image, set $\Sigma_s=\Sigma\times s$ 
and assume, in addition, that $\Sigma=\Sigma_0$ divides $U$ into two 
components: one bounded (i.e., such that its closure in $U$ is compact) 
and the other one, bounded or unbounded. 
(In many instances, this assumption can be omitted.)

The following, nearly trivial, proposition sums up a number of versions of the
nearby existence theorem:

\begin{Proposition}
\labell{prop:near-existence}
Assume that there exists an action selector $\sigma\colon \CH_W\to \R$
such that $\sigma(K)\leq C_U$ for every function $K\in \CH_U$ and
some constant $C_U$ independent of $K$. Then the nearby existence
theorem holds for proper functions $H$ on $U$: contractible in $W$ periodic 
orbits exist for a dense set of values of $H$. Equivalently, for every
shell $\{\Sigma_s\}$ there exists a value $s_0$ such that $\Sigma_{s_0}$ carries
a contractible in $W$ closed characteristic.\footnote{Depending on whether
this result is thought of in terms of $H$ or a shell, it is referred to as
a dense existence or nearby existence theorem.}
\end{Proposition}

As we will soon see, there exists an action selector on $\R^{2n}$ which
is \emph{a priori} bounded on bounded domains, \cite{HZ:book} and hence
the nearby existence theorem holds for $\R^{2n}$.

\begin{Remark}
In this proposition and in Proposition \ref{prop:almost-existence}, the
assumption that the shell divides $W$ is not essential; see \cite{SM}.
\end{Remark}

\begin{proof}[Proof of Proposition \ref{prop:near-existence}]
Given a shell $\{\Sigma_s\}$, consider a smooth function $F$ on $W$ such that
\begin{itemize}
\item $F$ is equal to some positive constant $C>C_U$ on the bounded
component of the complement to the shell and is identically zero on the
unbounded component,

\item within the shell, $F$ depends only on $s$ and is a monotone decreasing 
function with range $[0,C]$ and such that $0$ and $C$ are
the only critical values of $F$.
\end{itemize}
Note that the only critical values of $F$ on $W$ are again $0$ and $C$.
Furthermore, $0<\sigma(F)<C$ by \eqref{eq:positivity}
and since $\sigma(F)\leq C_U<C$. Hence, any contractible one-periodic 
orbit of $F$ with action $\sigma(F)$ must be 
non-trivial and lie on one of the regular levels of $F$, i.e., on some
hypersurface $\Sigma_{s_0}$. Hence, at least one of these hypersurfaces 
carries a contractible in $W$ closed characteristic.
\end{proof}

\subsubsection{The existence of an action selector}
\labell{sec:sel-exist}
Let us turn now to the problem of existence of an action selector. There are
numerous constructions of action selectors for different classes of 
manifolds $W$: in \cite{HZ:cap,HZ:book} an action
selector has been constructed for $W=\R^{2n}$ using a direct variational 
method on the space of $H^{1/2}$-loops; in \cite{Vi:functions} an action
selector has been defined for cotangent bundles using generating functions;
this construction has been extended to Lagrangian submanifolds in 
\cite{Oh1,Oh2} using Floer homology;
in \cite{Sc,Oh2} an action selector has been introduced for symplectically 
aspherical closed manifolds $W$, again by utilizing the Floer homology 
methods; these results have been extended to convex at infinity manifolds
in \cite{FS1}. Finally, in \cite{Oh3}, an action
selector was defined for closed symplectic manifolds which are not
necessarily convex or symplectically aspherical. 

Let us now outline the construction of an action selector, following 
mainly \cite{FS1,Sc}, for open symplectically aspherical manifolds 
which are convex at infinity.

Let $W$ be such a manifold. Recall that the Floer homology 
$\HF(H)$ for $H\in C^\infty_c(S^1\times W)$ is defined and independent of 
$H$, i.e., for any two 
functions these groups are canonically isomorphic, \cite{Fl1,Fl2}, see
also \cite{PSS,Sa,Sc:book}. Since $W$ is not closed but convex at
infinity, one has to extend $H$ to a function with a suitable growth at
infinity, without creating new periodic orbits, and then work with the
Floer homology of the resulting function; see \cite{FS1,Oh1,Vi:functors}.
Moreover, $\HF(H)$ is isomorphic to $H_*(W)$ up to a shift of degrees,
but this fact and the specific nature of the isomorphism are not very 
essential for us until later stages of the argument. At the moment,
one may interpret $H_*(W)$ as a common notation for the groups $\HF(H)$
identified with each other for different functions $H$. Furthermore, the
filtered Floer homology $\HF^{[a,\,b)}(H)$ is defined for any interval
$[a,\,b)$ with end-points outside of $\CS(H)$. This homology is constant
under deformations of $a$, $b$ and $H$ as long as $a$ and $b$ are outside
of the action spectrum of $H$, see, e.g., \cite{Vi:functors}.

There is a natural map 
$$
j^a\colon H_*(W)=\HF(H)\to \HF^{[a,\,\infty)}(H)
$$
induced by the quotient of complexes. Fix a non-zero element $u\in H_*(W)$
and set
$$
\sigma_u(H)=\inf\{a\mid j^a u=0\}.
$$
It is easy to see by using the invariance and monotonicity properties
of Floer homology that $\sigma_u$ meets the requirements (AS1) and (AS2).
Let $[\max]\in H_*(W)$ be the class of the maximum of a $C^2$-small
bump function. By calculating the Floer homology 
of such a bump function,
one can show that $[\max]\neq 0$. The selector $\sigma=\sigma_{[\max]}$
obviously satisfies (AS1)-(AS3). However, the $C^0$-continuity of $\sigma$,
or in general of $\sigma_u$, requires a proof (see, e.g., \cite{Oh2,Sc})
and this proof relies on the explicit construction of Floer's continuation
map.

The key to the proof is the following observation. Let $H_0\geq H_1$ be two
Hamiltonians whose one-periodic orbits are non-degenerate. Consider periodic
orbits $x_0$ of $H_0$ and $x_1$ of $H_1$ such that there exists a connecting
homotopy trajectory joining $x_0$ and $x_1$ for a certain ``linear'' 
monotone homotopy from $H_0$ to $H_1$. 
Then
\begin{equation}
\labell{eq:x0-to-x1}
A_{H_0}(x_0)- A_{H_1}(x_1)
\leq\int_0^1\max_{p\in W }\big(H_0(t,p)-H_1(t,p)\big)\,dt.
\end{equation}
This inequality can be obtained by a direct calculation and 
the $C^0$-continuity of $\sigma$ readily follows from \eqref{eq:x0-to-x1};
see \cite{Sc} for details. 

In a similar vein, \eqref{eq:x0-to-x1} implies the inequality
\begin{equation}
\labell{eq:max-to-sigma}
\sigma(H)\leq \int_0^1\max_{p\in W }H(t,p)\,dt,
\end{equation}
which can be established by setting $H=H_0$ in \eqref{eq:x0-to-x1},
taking a $C^2$-small function as $H_1$ and then letting $H_1$ go to zero.

\subsubsection{An \emph{a priori} bound for $\sigma$} 
\labell{sec:aprioribound}
One class of domains
$U$ for which $\sigma$ is \emph{a priori} bounded is the class of displaceable
domains, i.e., $U$ such that there exists a Hamiltonian symplectomorphism
$\varphi_H$ with $H\in \CH_W$ moving $U$ away from itself: 
$\varphi_H(U)\cap U=\emptyset$.
For example, bounded open subsets of $\R^{2n}$ are displaceable. 

To show that $\sigma$ is \emph{a priori} bounded on a displaceable domain
$U$, one argues as follows. 

Let $\varphi_H^t$ and $\varphi_K^t$ be time-dependent Hamiltonian flows
generated by $H$ and $K$ in $\CH_W$. Denote by $H\# K \in \CH_W$ 
the Hamiltonian generating $\varphi_H^t\circ\varphi_K^t$, i.e.,
$$
H\# K (t,p)=H(t,p)+K(t,(\varphi_H^t)^{-1}(p)).
$$
The crucial feature of the selector $\sigma$ defined in Section 
\ref{sec:sel-exist} is that $\sigma$ is sub-additive:
\begin{equation}
\labell{eq:subad}
\sigma(H\# K)\leq \sigma (H)+\sigma(K)
\end{equation}
or more generally $\sigma_{u\cap v}(H\# K)\leq \sigma_u (H)+\sigma_v (K)$,
where $u\cap v$ is the intersection of homology classes $u$ and $v$.

The sub-additivity of $\sigma$ is proved in \cite{Oh2,Sc} by making use of the 
pair-of-pants product introduced in \cite{Sc:thesis}; see also
\cite{PSS}. This is a product
$$
\HF^{[a,\,\infty)}(H)\otimes \HF^{[b,\,\infty)}(K)
\to \HF^{[a+b,\,\infty)}(H\#K)
$$
which is intersection of cycles on the full Floer homology $H_*(W)$ 
and such that the corresponding diagram commutes. 
(At this stage one specifically utilizes the isomorphism
$\HF(H)\to H_*(W)$ defined in \cite{PSS}, for this isomorphism sends the
pair-of-pants product to the intersection of cycles. However, 
it might also be possible to prove the sub-additivity
of $\sigma_{[\max]}$ by using the Floer continuation map.) Once the
existence of such a product is established, \eqref{eq:subad} follows from
the definition of $\sigma$; see \cite[Section 4]{Sc}.

The sub-additivity of $\sigma$ implies an \emph{a priori}
bound for $\sigma$. Let $K$ be a Hamiltonian whose time-one flow 
displaces $\supp H$, where, by definition, $\supp H=\bigcup_{t\in S^1}\supp H_t$. Denote by $K^-$ the Hamiltonian generating the flow $(\varphi_K^t)^{-1}$
so that $K\#K^-=0$. Then
\begin{equation}
\labell{eq:bound}
\sigma(H)\leq \sigma(K)+\sigma(K^-).
\end{equation}
To prove \eqref{eq:bound}, we argue as in \cite{FS1}. First observe that 
since $\varphi_K$ displaces $\supp H$, one-periodic orbits of $H\#K$ are 
exactly one-periodic orbits of $K$, and, as a consequence, 
$\CS(H\#K)=\CS(K)$. The same is, of course, true when $H$ is
replaced by the Hamiltonian $sH$ with $s\in \R$, i.e., $\CS((sH)\#K)=\CS(K)$.
By continuity of $\sigma$ and ``discontinuity'' of $\CS(K)$, we conclude that
$\sigma((sH)\# K)\in \CS((sH)\#K)=\CS(K)$ is independent of $s$. Setting
$s=0$ and $s=1$ yields $\sigma(K)=\sigma(H\#K)$. Therefore,
$$
\sigma(H)=\sigma(H\#K\#K^-)\leq \sigma(H\#K)+\sigma(K^-)
=\sigma(K)+\sigma(K^-),
$$
which proves \eqref{eq:bound}.

Among other important examples of displaceable domains $U$ are small 
neighborhoods of closed non-Lagrangian submanifolds of middle
dimension whose normal bundles
have non-vanishing sections, \cite{LS,pol:moving}. This (combined with
Macarini's stabilization trick, \cite{Ma}), leads to a local version 
of the nearby 
existence theorem for twisted cotangent bundles over surfaces (other 
than the torus) or, in higher dimensions, for exact magnetic fields;
see \cite{FS1}.

\subsection{The action selector method: closed manifolds}
\labell{sec:as-closed}
In this section we will briefly point out changes required in the 
action selector proof of the nearby existence theorem when $W$ is
a closed manifold.

First, let us observe that the construction of an action selector
outlined in Section \ref{sec:sel-exist} carries over to Hamiltonians
on closed manifolds $W$. (Here, as everywhere in this section, $W$ is
assumed to be symplectically aspherical.) Such an action selector has
the property $\sigma(H+\const)=\sigma(H)+\const$ and hence is never
\emph{a priori} bounded.

To circumvent this problem, one chooses a suitable normalization 
of Hamiltonians. (This is also necessary to ensure that $\sigma(H)$ 
depends only 
on $\varphi_H$.) However, in contrast with open manifolds where 
compactly supported Hamiltonians are automatically normalized, there appears to
be no natural choice of normalization for closed manifolds.

One possible choice is to restrict the action selector to Hamiltonians
vanishing on a neighborhood of a fixed point. This should allow one to
extend word-for-word the results and constructions of Section 
\ref{sec:as-open} to closed manifolds. However, certain details of
this approach are still to be worked out and we leave its discussion to
a later occasion. Here, we use the traditional normalization requiring
the mean value of Hamiltonians to be zero. More specifically,
let $\CH_W$ be the class of smooth functions on $S^1\times W$ 
normalized so that $\int_W H_t\omega^n =0$ and let $\CH_U$ be formed by 
Hamiltonians
$H\in  \CH_W$ such that $\supp X_{H_t}\subset U$ for all $t\in S^1$.

Proposition \ref{prop:near-existence} still holds when the class $\CH_U$
is defined in such a fashion. However, when $W$ is closed, obtaining an
\emph{a priori} bound for a selector on $\CH_U$ is more difficult than in 
the case of an open manifold, even though this still might be possible. 
We refer the reader to \cite{FGSS} for a detailed analysis of this question.
(The proof
of the \emph{a priori} bound from Section \ref{sec:aprioribound} does not
go through because inequality \eqref{eq:bound} need not hold. The reason
is that, even when $W\ssminus U$ is connected, $\CS(H\#K)$ will differ from
$\CS(K)$ by the constant equal to $H|_{W\ssminus U}$.)

Here, following \cite{Sc}, we choose a somewhat different approach.
Let, as in Section \ref{sec:sel-exist}, $\sigma_{[\max]}$ be the action
selector associated with the Floer homology class $[\max]\in H_*(W)$ 
corresponding to the maximum of a $C^2$-small bump function (shifted to have
zero mean). Set
$$
\gamma(H)=\sigma_{[\max]}(H)+\sigma_{[\max]}(H^-).
$$
Equivalently, $\gamma(H)=\sigma_{[\max]}(H)-\sigma_{[\min]}(H)$,
where $\sigma_{[\min]}$ is the action
selector associated with the Floer homology class $[\min]\in H_*(W)$ 
corresponding to the minimum of a negative $C^2$-small shifted bump function.
Similarly to action selectors for normalized functions, $\gamma(H)$ depends 
only on $\varphi_H\in \CalD$.

When $H$ is $C^2$-small, we have $\gamma(H)=\int_0^1(\max H-\min H)\,dt$. 
Furthermore, 
$\gamma$ is still $C^0$-continuous in $H$, sub-additive, but not necessarily 
monotone. Then, since $H_t^-(p)=-H_t\left((\varphi_H^t)^{-1}(p)\right)$, 
inequality \eqref{eq:max-to-sigma} translates to
\begin{equation}
\labell{eq:var-to-gamma}
\gamma(H)\leq \int_0^1\big(\max_{p\in W }H(t,p)-\min_{p\in W }H(t,p)\big)\,dt.
\end{equation}
Note also that $\gamma(K^-)=\gamma(K)$. As we have pointed out above,
the proof of \eqref{eq:bound} given in Section \ref{sec:aprioribound} does 
not go through for the normalization we use here. However, a similar but 
more involved argument proves the upper bound
\begin{equation}
\labell{eq:bound-gamma}
\gamma(H)\leq 2\gamma(K),
\end{equation}
where the time-one flow of $K$ displaces $U \supset \supp(X_H)$.
(Here, by definition, $\supp(X_H)=\bigcup_{t\in S^1}\supp X_{H_t}$.)
We refer the reader to \cite{Sc} for detailed proofs of these facts.

Observe now that Proposition \ref{prop:near-existence} holds when 
an \emph{a priori} bounded selector $\sigma$ is replaced by 
$\gamma$ which satisfies \eqref{eq:var-to-gamma} and is \emph{a priori}
bounded on $U$. (To see this, only minor modifications in the proof
are needed.) As a consequence of \eqref{eq:bound-gamma}, we immediately 
obtain the nearby existence 
theorem for displaceable domains in compact symplectically aspherical 
manifolds.

\subsection{Limitations of the action selector method}
The class of manifolds $W$ and domains $U$ to which the action
selector method applies is rather limited. For example, a bounded selector 
or a bounded function $\gamma$ do not exist for $U=W=\T^2$ 
(nor for $\T^4$ or 
$\T^2\times S^2$ with split symplectic structure). Indeed, it is
easy to see that there is a function $H\in\CH_{\T^2}$ with
$\CS(H)=\{\min H,\max H\}$ and such that $\sigma_{[\max]}(H)=\max H$ and
$\sigma_{[\min]}(H)=\max H$ are both arbitrarily large. 
The same is true for tubular neighborhoods of the zero section in
cotangent bundles of some compact manifolds (e.g., of surfaces with genus
$g\geq 1$). As a consequence, 
Proposition \ref{prop:near-existence} cannot be applied to prove
nearby existence on these manifolds as long as the standard Floer homology 
(even accounting for non-contractible orbits) is employed. This is 
the case already for $W=\T^2$ even though every regular level of $H$ is  
comprised of (not-necessarily 
contractible) periodic orbits. The problem is that the periodic orbit can
``migrate'' from one homotopy class to another depending on the function
in question. 

Another class of manifolds to which this method does not
apply is formed by geometrically bounded, but not convex manifolds. Among such
manifolds are many twisted cotangent bundles and also universal coverings
of some manifolds, which justifies the interest in this class.

The problem is that there is no known satisfactory definition of an action 
selector for general geometrically bounded manifolds. For example, the 
obstacle in the homological 
approach is that $\HF^{[a,\, b)}(H)$ for $H\in \CH_W$
has been defined only for intervals $[a,\,b]$ 
which do not 
contain zero, while the homological definition of $\sigma$ requires 
$\HF^{[a,\, b)}(H)$ to be defined for all intervals. One possible solution 
is to consider an action selector only on the class of non-negative 
(autonomous or time-dependent) Hamiltonians or a yet more narrow 
class of functions (containing the functions $F$ from the proof of 
Proposition \ref{prop:near-existence}) and to extend the homological 
definition to this class. The proof of the
\emph{a priori} bound outlined in this section will not carry over to such a 
narrow class, but a different argument (e.g., akin to the symplectic homology
method from Section \ref{sec:symhom}) may. 

Furthermore, when a geometrically bounded manifold is ``asymptotically
convex'', as are for instance twisted cotangent bundles, one should
be able to define Floer homology for all intervals of action by suitably
extending the function at infinity as in the convex case. When this is done,
many of the results that hold for convex manifolds should remain valid, 
although certain technical difficulties have to be dealt with. 

However, as of today, no work in either of these directions has been carried
out and it is not clear whether or not this would lead to new results.

\subsection{Symplectic homology} 
\labell{sec:symhom}
A slightly different approach to proving the nearby existence theorem, 
which also utilizes Floer homology, 
relies on symplectic homology introduced and investigated in 
\cite{CFH,CFHW,FH,FHW}.

Let, as above, $U$ be a domain in a symplectic manifold $W$. The
symplectic homology of $U$ is defined as 
$$
\SH^{[a,\,b)}(U)=\varprojlim_{H}\HF^{[a,\,b)}(H),
$$
where the inverse limit is taken over all $H\in C^\infty_c(S^1\times U)$. 

The symplectic homology of $U$ detects closed characteristics in an
arbitrarily narrow thickening $\{\Sigma_s\}$ of $\Sigma=\partial U$,
provided that $\Sigma$ is smooth. In particular, when
$\SH^{[a,\,b)}(U)\neq 0$, in any thickening $\{\Sigma_s\}$ there exists a
hypersurface $\Sigma_{s_0}$ carrying a closed characteristic. (This
fact is an essentially immediate consequence of the definition; see
\cite{FH}.)  Therefore, the nearby existence theorem holds for a
function $H$ whenever the symplectic homology of the sublevels $\{H<c\}$
is non-trivial. The problem is thus reduced to verifying that 
symplectic homology does not vanish.

Let us show how this approach works,  for example,  when $W=\R^{2n}$.
Let $U$ be a bounded domain in $\R^{2n}$ which we assume to contain
the origin. Thus, there exist two open balls $B_r$ and $B_R$ centered
at the origin and such that $B_r\subset U\subset B_R$.
The symplectic homology has a natural monotonicity property: an inclusion
$U\subset V$ induces a map $\SH^{[a,\,b)}(V)\to \SH^{[a,\,b)}(U)$. The map
$\Psi\colon \SH^{[a,\,b)}(B_R)\to \SH^{[a,\,b)}(B_r)$ factors as
$$
\SH^{[a,\,b)}(B_R)\to \SH^{[a,\,b)}(U)\to \SH^{[a,\,b)}(B_r),
$$
and hence $\SH^{[a,\,b)}(U)\neq 0$ when $\Psi\neq 0$. The symplectic homology
of a ball and the map $\Psi$ can be calculated explicitly (see \cite{FHW,Her1})
and, indeed, it turns out that $\Psi\neq 0$ for some $b>a>0$ and 
the nearby existence 
in $\R^{2n}$ follows. Today, this calculation can be carried out particularly 
easily (cf. \cite{BPS,CGK}) if one makes use of Po\'zniak's theorem giving
Floer homology in the Morse--Bott non-degenerate case, \cite{poz}.

This method also applies in the setting of the generalized Weinstein--Moser
theorem discussed in the introduction. Let $U$ be 
a neighborhood of a closed symplectic submanifold $M$ of a geometrically 
bounded symplectically aspherical manifold $W$. Then, by taking suitably 
defined symplectic tubular neighborhoods of $M$ as $B_r$ and $B_R$, one
can show that $\SH^{[a,\,b)}(U)\neq 0$, \cite{CGK}. This leads to a proof
of a nearby existence theorem for narrow shells enclosing $M$ or small
values of functions having a minimum (say, equal to zero) along $M$. 
Note that, since here we utilize Floer homology only for intervals that
do not contain zero, the method can be used for geometrically bounded
manifolds which are not necessarily convex.

Furthermore, this method can also be cast in the framework of energy selectors
defined on a suitable class of functions (cf. \cite{Her1}). However, 
the benefits of this approach are unclear and we omit the details here.

\begin{Remark}[Stability of the area spectrum, \cite{CFHW}]
As a side remark, let us mention one more application of symplectic
homology and Po\'zniak's theorem. Here, again, we focus on the main
idea of the argument rather than on a complete proof.
Let $\Sigma$ be a smooth hypersurface in $W$.
The area (or action) spectrum $\CA(\Sigma)$ is the collection of
symplectic areas bounded by contractible closed characteristics on $\Sigma$,
including iterated closed characteristics. It is not hard to see using
Po\'zniak's theorem that, when $\Sigma$ is the boundary of an open domain $U$,

$$
\CA(\Sigma)=\{a\in\R\mid \SH^{[a-\eps,\,a+\eps)}(U)\neq 0\text{ for any
small $\eps>0$}\},
$$
provided that $\Sigma$ has contact type and all (iterated)
closed characteristics on $\Sigma$ are non-degenerate. 
(An additional, more subtle, argument is needed here when
different closed characteristics bound equal areas.)
Observing that the right hand side of this
equality depends solely on $U$, we conclude that under these hypotheses 
$\CA(\Sigma)$
is determined by $U$. This is the stability of the area spectrum theorem,
\cite{CFHW}. More specifically, let $\Sigma$ and $\Sigma'$ be smooth
contact type hypersurfaces in $W$ bounding open domains $U$ and $U'$, 
respectively. Assume that all (iterated) closed characteristics on $\Sigma$ 
and $\Sigma'$ are non-degenerate and that there exists a symplectomorphism
of $W$ sending $U$ to $U'$. Then, $\CA(\Sigma)=\CA(\Sigma')$.
\end{Remark}

\subsection{Other applications} In this section we briefly discuss
some other applications of the action selector, which are important for
what follows.

\subsubsection{Homological capacity}
\labell{sec:hom-cap}
An action selector $\sigma$ can be used to introduce an 
invariant, the homological capacity $\wsh(U)$, of a domain $U$ in a symplectic 
manifold $W$. Assume first that $W$ is open and convex at infinity. To an
action selector $\sigma$ on $W$, we associate a function $\wsh$ on
open subsets of $W$ by setting
$$
\wsh(U)=\sup\{\sigma(H)\mid H\in \CH_U\}\in (0,\infty].
$$
By definition, $\wsh(U)<\infty$ if and only if $\sigma$ is \emph{a priori}
bounded on $U$, and hence the nearby existence theorem holds for $U$.

When $W$ is closed, we set
$$
\wsh(U)=\sup\{\gamma(H)\mid H\in \CH_U\}\in (0,\infty].
$$
The homological capacity is 
invariant under symplectomorphisms of $W$,
monotone with respect to inclusions of domains, and homogeneous of
degree one with respect to scaling of $\omega$, i.e.,
$$
\wsh(U,\lambda\omega) = |\lambda|\wsh(U,\omega)
\text{ for any non-zero $\lambda\in \R$.}
$$ 
Here we do not touch upon the general definition and properties of symplectic
capacities and refer the reader to, e.g., 
\cite{Ha,Her1,Her2,Ho:cap,HZ:cap,HZ:book}
for a detailed discussion of this subject.

In what follows we will always assume that $\wsh$ is associated to
$\sigma=\sigma_{[\max]}$ or $\gamma=\sigma_{[\max]}-\sigma_{[\min]}$.
Then, it is not hard to show that $\wsh(B_R)=\pi R^2$ by 
calculating the Floer homology of bump functions on $B_R$ (see \cite{FHW}
and Section \ref{sec:symhom}).
Similarly, $\wsh(S^2)=\area(S^2)$. However, $\wsh(W)=\infty$, when $W$
is a closed orientable surface other than $S^2$ or $W=\T^2\times S^2$ or
$W=\T^4$.

In Section \ref{sec:comparison} we will utilize the homological capacity 
$\wsh$ as an upper bound for the Hofer--Zehnder capacity.

\subsubsection{Hofer's geometry} Denote by $\CalD$ the universal cover
of the group of Hamiltonian symplectomorphisms of $W$. To be more precise,
$\CalD$ is  formed by time-one flows $\varphi_K$ of 
Hamiltonians $K\in \CH_W$, where each $\varphi_K$ comes 
together with the homotopy class (with fixed end points) of the path
$\varphi^t_H$ for $t\in [0,\,1]$. 
Hofer's norm $\parallel \cdot \parallel_{\H}$ is defined as
$$
\parallel K \parallel_{\H}=\int_0^1 (\max K_t-\min K_t)\, dt
$$
for $K\in \CH_W$. For $\psi\in\CalD$, set
$$
\rho(\psi)=\inf\{\parallel K \parallel_{\H}\,\mid \text{ $K$ generates
$\psi$, i.e., $\psi=\varphi_K$}\}.
$$
It is easy to see that $\rho(\psi,\varphi):=\rho(\psi\varphi^{-1})$ is
a bi-invariant metric on $\CalD$, provided that $\rho$ is 
non-degenerate, i.e.,
$$
\rho(\psi)>0\quad\text{ iff }\quad \psi\neq\id .
$$
It turns out that $\rho$, known as Hofer's metric, is indeed non-degenerate 
for any symplectic manifold $W$. Non-degeneracy of
$\rho$ has been established through a series of more and more general results
starting with $W=\R^{2n}$ (see e.g., \cite{HZ:book}) and ending with
the proof for an arbitrary $W$ in \cite{LMD}. We refer the reader to 
\cite{polt:book} for an introduction to Hofer's geometry and further 
references. Here we only outline the 
proof of non-degeneracy for symplectically aspherical convex manifolds.

Assume first that $W$ is open and convex at infinity.
Let $\varphi_K\neq \id$. Pick an open set $U$ displaced by
$\varphi_K$, i.e., such that $\varphi_K(U)\cap U=\emptyset$, and let
$H\in\CH_U$. The idea is to again utilize \eqref{eq:bound}, but this time
to obtain a lower bound for the right hand side. Namely, 
\eqref{eq:max-to-sigma} yields the inequality 
$\sigma(K)\leq \int_0^1\max K_t\,dt$ and also
$\sigma(K^-)\leq -\int_0^1\min K_t\,dt$, as can be
easily seen from the definition of $K^-$. Therefore, 
by \eqref{eq:bound},
$$
\sigma(H)\leq \parallel K\parallel_{\H},
$$
and hence, once $\supp H \subset U$ and $\sigma(H)>0$,
\begin{equation}
\labell{eq:nondeg}
0<\sigma(H)\leq \rho(\varphi_K).
\end{equation}
This proves non-degeneracy of $\rho$ for convex at infinity 
symplectically aspherical open manifolds. Furthermore, let us define
the displacement energy of $U$ as $e(U)=\inf\rho(\psi)$, where 
$\psi\in \CalD$ displaces $U$. Then \eqref{eq:nondeg} translates into 
the upper bound
\begin{equation}
\labell{eq:displ}
\wsh(U)\leq e(U),
\end{equation}
which is a minor improvement over the upper bound
$\wsh(U)\leq 2e(U)$ established in \cite{FS1}.

When $W$ is closed and symplectically aspherical, we argue in a similar 
fashion. Let, as above, $\varphi_K$ displace $U$. Note that
$\gamma(K)\leq \rho(\varphi_K)$, by \eqref{eq:var-to-gamma}. Hence, we 
conclude from \eqref{eq:bound-gamma} that
$$
\gamma(H)\leq 2\rho(\varphi_K)
$$
for any $H\in\CH_U$. When $H$ is a $C^2$-small 
nonzero function, $\gamma(H)=\int_0^1 (\max H_t-\min H_t)\,dt >0$. Therefore,
$\rho(\varphi_K)>0$, which proves non-degeneracy of $\rho$. Furthermore, we
also obtain the upper bound, \cite{Sc}:
\begin{equation}
\labell{eq:displ2}
\wsh(U)\leq 2e(U).
\end{equation}
In fact, the more accurate upper bound \eqref{eq:displ} still holds in this 
case; see \cite{FGSS}.

Non-degeneracy of $\rho$ for an arbitrary
symplectic manifold has been established in \cite{LMD} by showing, using 
entirely different methods, that
$$
\CGR(U)\leq 2\rho(\varphi_K),
$$
where $\CGR(U)$ is the Gromov capacity of $U$, i.e.,
$$
\CGR(U)=\sup \{\pi R^2\mid\text{$B_R$ symplectically embeds into $U$}\}.
$$
Note also that since $\CGR(U)\leq \wsh(U)$, \eqref{eq:nondeg}
yields that in fact $\CGR(U)\leq \rho(\varphi_K)$, whenever $W$ is 
symplectically aspherical, open, and convex at infinity.

\section{The almost existence theorems}
\labell{sec:almost}
Virtually all known proofs of the almost existence theorems are based
on the notion of the Hofer--Zehnder capacity, \cite{HZ:book}, with the
exception of \cite{St} which precedes this notion.

\subsection{The Hofer--Zehnder capacity and almost existence}
Let $V$ be a symplectic manifold without boundary. Denote by $\CHHZ(V)$ the 
class of smooth non-negative functions $K$ on $V$ such that
\begin{itemize}
\item $K$ is compactly supported if $V$ is not closed or $K$ vanishes on
some open set if $V$ is compact,
\item $K$ is constant near its maximum.
\end{itemize}
We will refer to such functions as Hofer--Zehnder functions. Also,
let us call a non-trivial periodic orbit of $K$ with period $T\leq 1$ fast.
Otherwise, an orbit will be called slow. A Hofer--Zehnder function without 
non-trivial fast periodic orbits will be said to be admissible.
Following \cite{HZ:cap,HZ:book}, recall that the Hofer--Zehnder capacity 
of $V$ is 
defined as 
$$
\CHZ(V)=\sup\{\max K \mid K\in \CHHZ(V)\text{ and $K$ 
is admissible}\} 
\in (0,\infty].
$$
The capacity $\CHZ(V)$ does not change when the
assumptions that $K$ is non-negative and/or that $K$ is constant near
its maximum are dropped, \cite{gg3}. The Hofer--Zehnder capacity has the
same general properties as the homological capacity.
One should think of $\CHZ$ as a higher dimensional analogue of the area.
We will elaborate on this point in Section \ref{sec:cap-fun} and here
we only mention that by
using the area--period relation (see Section \ref{sec:cap-fun}) one can
show that for closed orientable surfaces $\CHZ$ is exactly equal to the 
area, \cite{Si}, in contrast with $\wsh$.

The following result asserts that to prove almost existence in $V$,
it suffices to establish that $\CHZ(V)$ is finite.

\begin{Proposition}[\cite{HZ:book}]
\labell{prop:almost-existence}
Assume that $\CHZ(V)<\infty$. Then the almost existence theorem holds for 
proper $C^2$-functions $H$ on $V$: periodic orbits of $H$ 
exist on almost all, in the sense of measure theory, regular levels of $H$. 
Equivalently, for every shell $\{\Sigma_s\}$ which bounds a domain in $V$, the 
hypersurfaces $\Sigma_{s}$ carry closed characteristics for almost all $s$.
\end{Proposition}

This proposition is a rather simple consequence of the definition of $\CHZ$ and
the Arzela--Ascoli theorem; see \cite{HZ:book}. The assumption that
$\{\Sigma_s\}$ bounds a domain is superfluous, \cite{SM}.
In Section \ref{sec:cap-fun} we will prove a more precise version of 
Proposition \ref{prop:almost-existence}.

One can also incorporate the homotopy class of an orbit in the almost existence
theorem and in the definition of $\CHZ$, \cite{Sc}. The simplest way to do 
this is as follows. Fix an ambient symplectic manifold $W$. Let us now
modify the definition of the Hofer--Zehnder capacity by requiring $V$ to
be an open subset of $W$ and requiring $K\in \CHHZ(V)$ to have no 
non-trivial contractible in $W$ fast periodic orbits.
We denote the resulting capacity by $\CW(V)$. Then Proposition
\ref{prop:almost-existence} holds for contractible in $W$ periodic orbits
provided that $\CW(V)<\infty$. Clearly, $\CHZ(V)= \CW(V)$, when
$W$ is simply connected, and 
\begin{equation}
\labell{eq:cap-contr}
\CHZ(V)\leq \CW(V)
\end{equation}
in general. The strict inequality is possible: for example, 
$\CHZ(S^1\times (0,1))=\area(S^1\times (0,1))<\CW(S^1\times (0,1)) =\infty$,
where we view the annulus $S^1\times (0,1)$ as a subset of itself. This
also shows that the choice of the ambient manifold, even though it is not
included in the notation, effects the value of $\CW$. (Replace
$W=S^1\times (0,1)$ by $W=\R^2$.) In this connection, we also note that
$\CHZ(U)<\infty$ for any bounded open subset $U$ of $T^*\T^n$, \cite{ji},
while $\CW(U)=\infty$ when $U$ contains the zero section.
Indeed, $U$ can be symplectically embedded into a bounded subset of $\R^{2n}$
since $\T^n$ admits a Lagrangian embedding in $\R^{2n}$. In particular, the
almost existence theorem holds in  $T^*\T^n$.

Beyond dimension two, little is known about the capacity $\CHZ$ in contrast
with the capacity $\CW$. For example,
all results discussed in what follows deal with contractible periodic 
orbits and thus concern the capacity $\CW$. Of course, by
\eqref{eq:cap-contr}, an upper bound for $\CW$ implies an upper bound for
$\CHZ$.

\subsection{Finiteness of the Hofer--Zehnder capacity}

In view of Proposition \ref{prop:almost-existence}, to prove the almost
existence theorem in $V$ it suffices to show that $\CHZ(V)<\infty$. Let us
start by proving that the Hofer--Zehnder capacity of a bounded domain in
$\R^{2n}$ is finite. 

\begin{Theorem}[\cite{HZ:cap,HZ:book}]
\labell{thm:capacity}
Let $U$ be a bounded domain in $\R^{2n}$. Then $\CHZ(U)<\infty$. Moreover,
$\CHZ(B_R)=\pi R^2$, where $B_R$ is the ball of radius $R$.
\end{Theorem}

\begin{proof}
By monotonicity, the first assertion is an immediate consequence of the 
inequality
$\CHZ(B_R)<\infty$. To show that $\CHZ(B_R)=\pi R^2$, it suffices to prove
that $\CHZ(B_R)\leq \pi R^2$; the opposite inequality is easy. (It is usually
the case that establishing an upper bound for the capacity is hard while a 
lower bound can be easily obtained by definition.) The key is the following
result:

\begin{Proposition}[\cite{HZ:book}]
\labell{prop:capacity}
Let $H$ be a Hofer--Zehnder function on $B_R$ such that 
$\max H >\pi R^2$. Then the Hamiltonian
flow of $H$ has a non-trivial one-periodic orbit.
\end{Proposition}

Following \cite{gg3}, let us outline a proof of the proposition which
relies on the calculation of the Floer homology and generalizes to some
other situations.

Let $H$ be as in Proposition \ref{prop:capacity}. There exist non-negative
functions $K^-\leq H\leq K^+$, supported in $B_R$ and depending only on the
distance to the origin, and such that the action of $K^\pm$ on the nearest
to the origin non-trivial one-periodic orbits is greater than $\max K^+$.
(These orbits are Hopf circles on a sphere enclosing the origin.) To find
such functions, we make use of the assumption that $H$ is constant near 
its maximum and take, as $K^\pm$, functions that squeeze $H$ from above 
and below as tightly as possible and that depend  only on 
the distance to the origin. Now, applying Po\'zniak's theorem \cite{poz},
one can show that 
$$
\HF^{[a,\,\infty)}_{m}(K^\pm)=\Z_2
$$
for some $m$ and $a>\max K^+\geq \max H$. This homology group is ``generated''
by the non-trivial one-periodic orbits closest to the origin. 
Moreover, by examining the natural homotopy from $K^+$ to $K^-$, one can prove
that the monotonicity map, which factors as
$$
\Z_2=\HF^{[a,\,\infty)}_{m}(K^+)\to\HF^{[a,\,\infty)}_{m}(H)\to
 \HF^{[a,\,\infty)}_{m}(K^-)=\Z_2,
$$
is an isomorphism. Hence, $\HF^{[a,\,\infty)}_{m}(H)\neq 0$. Every trivial
periodic orbit of $H$ has action in the interval $[0,\max H]$. Thus, since
$a>\max H$, the flow of $H$ must have a non-trivial one-periodic orbit.
This concludes the proof of Proposition \ref{prop:capacity} and Theorem
\ref{thm:capacity}.
\end{proof}

Observe that Proposition \ref{prop:capacity} is stronger than what is
needed: the proposition
guarantees the existence of a non-trivial orbit with period $T=1$ while
$T\leq 1$ is sufficient  to prove Theorem \ref{thm:capacity}. Yet, this fact, 
somewhat surprisingly, appears to have no interesting applications.

Also note that there are compact aspherical symplectic manifolds with
infinite Hofer--Zehnder capacity. The basic example of such a manifold
is Zehnder's torus, \cite{Ze}, i.e., a torus $\T^{2n}$ with an
irrational symplectic structure. All other examples of such manifolds known
to the author are derivatives of Zehnder's torus. As has been pointed out
in the introduction, the nearby/almost existence theorem also fails for 
Zehnder's torus. 

Let us now elaborate on some general principles concerning the almost
existence or nearby existence theorem and the Weinstein conjecture. The
existence problem for periodic orbits on a given energy level is
dual in a certain, rather vague, sense to the existence problem for periodic 
orbits of a fixed period. The orbits of a fixed period are known to exist 
(Arnold's conjecture) and a variety of methods (e.g., Floer homology) 
of proving this fact have been developed during the last 
two decades. In contrast with this, unless the function is assumed to be convex
(and then the duality acquires a precise meaning), the problem for a fixed
energy level may fail to have a solution,  as counterexamples to
the Hamiltonian Seifert conjecture show. Moreover, there seems to be no direct
method to tackle the problem of existence of periodic orbits for a 
sufficiently large set of energy values (e.g., dense or full measure). All 
known methods, including those outlined above, first reduce the problem
to proving that the flow of a function $H$ with sufficiently large 
variation $\max H-\min H$ possesses (non-trivial, fast, or period one) 
periodic orbits. Thus, all proofs of these
theorems hinge on a general principle that a compactly supported function  
with sufficiently large variation must have fast non-trivial periodic orbits
or even one-periodic orbits if the function is constant near its maximum. This
principle, which is often true but fails in general (Zehnder's torus!), 
can already be proved for some manifolds by utilizing, for 
example, Floer homology.

Coming back to the discussion of the Hofer--Zehnder capacity, let us
mention just one more of its properties. Namely, $\CHZ(Z_R)=\pi R^2$, where
$Z_R=D^2_R\times \R^{2n-2}\subset \R^{2n}$ is a symplectic cylinder, 
\cite{HZ:cap,HZ:book}. (Thus, the Hofer--Zehnder capacity of $Z_R$ is 
finite even though its volume is infinite.) This can be easily proved 
as follows, 
\cite{gg3}. First note that our proof of Proposition \ref{prop:capacity} 
readily extends to ellipsoids and shows that the capacity of an ellipsoid is 
equal to its minimal cross-section area (by a plane through the origin). 
Then exhausting $Z_R$ by more and 
more elongated ellipsoids we see that $\CHZ(Z_R)=\pi R^2$.

The next application of our method concerns the generalized Weinstein--Moser
conjecture, the motion of a charge in a magnetic field, and the relative
almost existence theorem. Let $U$ be 
a neighborhood of a closed symplectic submanifold $M$ of a geometrically 
bounded symplectically aspherical manifold $W$. Then, we can emulate the
proof of Proposition \ref{prop:capacity} by taking, as $K^\pm$, functions
that depend only on the distance to $M$ in a suitable metric and that
are supported
near $M$. As a consequence, we obtain the relative almost existence
theorem:  near $M$, almost all levels of a function $H$ attaining a 
minimum along $M$ carry contractible periodic orbits, \cite{gg3}. 
(Note that this, in turn, requires replacing the ordinary Hofer--Zehnder 
capacity in Theorem \ref{thm:capacity} by its relative counterpart.
 The relative Hofer--Zehnder capacity $\CHZ(W,M)$, where
$M\subset W$ is a compact subset, is defined just as the ordinary capacity,
but the function $K$ is required to attain its maximum along $M$, see, e.g.,
\cite{gg3,HZ:cap,la} for more details.) 

Furthermore, a more refined version of the Floer homology calculation
outlined here 
can be used to show that $\CW(U)<\infty$ when $U$ is a small neighborhood
of a closed symplectic submanifold $M$ with homologically trivial normal
bundle in a closed symplectically aspherical manifold, \cite{Ke:new}. This
implies the almost existence theorem in $U$.

\subsection{Comparison of $\CW$ with $\wsh$} 
\labell{sec:comparison}
In this section we will outline a different approach to
obtaining an upper bound for $\CHZ$. To be more specific, we prove that
$\CW(U)\leq \wsh(U)$. In our proof of this inequality we draw heavily
on \cite{HZ:book,FS1,Sc}, but some details appear to be new.

Let $W$ be an open manifold (with $\omega|_{\pi_2(W)}=0$) on which a 
continuous selector $\sigma$ 
satisfying (AS1) and (AS2) is defined. Assume also that
\begin{equation}
\labell{eq:sel-small}
\sigma(H)=\max H,\text{ when $H$ is $C^2$-small and independent of time}.
\end{equation}
This requirement strengthens (AS3) and holds for
the energy selector defined in Section \ref{sec:sel-exist}
because the Floer complex of a $C^2$-small autonomous Hamiltonian 
is equal to its Morse complex; see \cite{FHS}.
Let $\wsh$ be the associated homological capacity. 

\begin{Proposition}[\cite{FS1,Sc}]
\labell{prop:cap-ineq}
$\CW(U)\leq \wsh(U)$ for any $U\subset W$.
\end{Proposition}

This, by \eqref{eq:cap-contr}, implies that $\CHZ(U)\leq \wsh(U)$.

\begin{proof} 
The proposition is an immediate consequence of the following

\begin{Lemma}
\labell{lemma:sigma}
Let $H\in \CH_W$ be an autonomous Hamiltonian whose flow has no fast 
non-trivial contractible periodic orbits. Then $\sigma(H)=\max H$.
\end{Lemma}

To derive the proposition from the lemma, we just use the definition
of the capacities:
$$
\wsh(U)=\sup\{\sigma(H)\mid H\in \CH_U\}
\geq \sup\{\max H\mid\text{$H$ as in the lemma}\}
=\CW(U).
$$

The assertion of the lemma is very plausible. When $W$ is convex at infinity
and $\sigma=\sigma_{[\max]}$, one can argue
as follows. By continuity, it suffices to prove the lemma for a generic 
$H$. Then, since $H$ does not have  non-trivial
contractible periodic orbits, the Floer complex of $H$ is generated by
the critical points of $H$ and hence $\sigma(H)$ is a critical value of $H$.
(Moreover, without loss of generality we
may assume that the eigenvalues of $d^2H$ at the critical points of $H$ 
are small. Indeed, this can be achieved by replacing $H$ by $f\circ H$,
where $f\colon \R\to \R$ is a suitable diffeomorphism which is $C^0$-close
to identity, see \cite{gg3,FS3}. Then $\sigma(H)$ must be the value of $H$
at one of its local maxima, for maxima are the only critical points of $H$
with correct index. Hence, $\max H=\sigma(H)$ if $H$ has a unique maximum.)
However, to show that $\sigma(H)$ is delivered by the global maximum,
one has to analyze the Floer differential of $H$ and this indeed can be done; 
see, e.g., \cite{FS1,Sc}. At this point one considers the homotopy $\lambda H$ 
with $\lambda\in (0,1]$ and fully uses the assumptions on $H$. Note also that
it seems to be unknown whether in the hypotheses of the lemma fast
orbits can be replaced by one-periodic orbits.

One can also prove the lemma in a formal set-theoretic way, using
only the continuity of $\sigma$, (AS1), (AS2) and \eqref{eq:sel-small}.
Denote by $\CS$ the action spectrum $\CS(H)$.
Observe that by the hypothesis of the lemma, every contractible one-periodic 
orbit of $\lambda H$ is trivial for all $\lambda\in (0,1]$, and hence 
the action spectrum $\CS(\lambda H)$ is again comprised entirely
of the critical values of $\lambda H$. Therefore, $\CS(\lambda H)=\lambda \CS$.
Set
$$
\sigma(\lambda)=\sigma(\lambda H)\in \lambda\CS.
$$
Then $\sigma(\lambda)/\lambda$ is a continuous function of $\lambda$ with
values in $\CS$. Since the action spectrum $\CS$ is a nowhere dense,
$\sigma(\lambda)/\lambda=\const$. For a small $\lambda>0$, the
function $\lambda H$ is $C^2$-small and hence 
$\sigma(\lambda H)=\max(\lambda H)=\lambda\max H$. It follows that
$\sigma(\lambda)=\lambda\max H$ for all $\lambda\in (0,\,1]$ which 
concludes the proof of the lemma and of the proposition. 
\end{proof}

Proposition \ref{prop:cap-ineq} also holds (see \cite{Sc}) when
$W$ is closed and symplectically aspherical and $\wsh$ is the homological
capacity defined in Section \ref{sec:hom-cap} for
$\gamma=\sigma_{[\max]}-\sigma_{[\min]}$. This can be proved essentially
in the same way as the result for open manifolds.

Combining Proposition \ref{prop:cap-ineq} with \eqref{eq:displ} and 
\eqref{eq:displ2}, we infer that a displaceable open set has a finite 
Hofer--Zehnder capacity:
\begin{eqnarray}
\labell{eq:cap-displace1}
\CW(U)&\leq & e(U)\text{ if $W$ is open,}\\
\labell{eq:cap-displace1a}
\CW(U)&\leq & 2e(U)\text{ if $W$ is closed,}
\end{eqnarray}
whenever $W$ is convex and symplectically aspherical. 
These are versions of the so-called ``capacity--displacement energy
inequality''. In particular, the almost existence theorem holds for 
displaceable sets. For example,
let $L\subset W$ be a  closed non-Lagrangian submanifold of middle
dimension whose normal bundle has a non-vanishing section. Assume also that
on $W$ there exists a selector $\sigma$ satisfying \eqref{eq:sel-small}. Then, 
by \cite{LS,pol:moving},
a small neighborhood $U$ of $L$ is displaceable and hence
$\CHZ(U)<\infty$. As in Section \ref{sec:aprioribound},
this fact leads to a local version of the almost 
existence theorem for twisted cotangent bundles over surfaces (other 
than the torus) or, in higher dimensions, for exact magnetic fields;
see \cite{FS1, Ma}.

Proposition \ref{prop:cap-ineq} explains why 
the almost existence theorem can often be proved in the same setting as
that of the weaker nearby existence theorem. Indeed, most of the proofs of the
former, with some notable exceptions, rely on first establishing that
a version of the homological capacity is finite. However, in this case
the Hofer--Zehnder capacity is also finite, which implies almost existence 
theorem.

\begin{Remark}
Proposition \ref{prop:capacity} suggests that in the definition
of the Hofer--Zehnder capacity one may replace the condition that
$K$ has no fast periodic orbits by a weaker requirement that it has no
one-periodic orbits (cf. \cite{Ke:new}). The resulting invariant is still 
a capacity. For this capacity
the upper bound from Proposition \ref{prop:cap-ineq} may fail to hold,
unless an additional requirement on the critical set of $K$ is imposed,
cf. \cite{FGSS}. However, this capacity is \emph{a priori}
bounded on $\R^{2n}$ and sufficient for the proof of the almost existence
theorem. Furthermore, the inequality $T\leq \LD \cf_H(h)$ in
Proposition \ref{prop:apr2} is then replaced by the equality. (In fact, 
when $\CHZ$ is modified in this way, for
every converging sequence 
$\big(\cf_H(h)-\cf_H(h_i)\big)/(h-h_i)\to A$ there exists
an $A$-periodic orbit (not necessarily simple) on the level $H=h$.) As
a consequence, this capacity is represented on hypersurfaces of
restricted contact type in $\R^{2n}$ (cf. Proposition \ref{prop:repr}).
\end{Remark}

\subsection{Holomorphic curves and almost existence}
The holomorphic curve approach to the nearby existence and almost existence
theorems goes back to \cite{FHV,HV2}. Here we outline only the
main idea of the method, following essentially \cite{HV2} and
omitting technical details which are in some instances quite involved. 

Let $(W,\omega)$ be a closed symplectic manifold and let $L_-$ and $L_+$ be 
two disjoint
closed submanifolds of $W$. Fix a generic almost complex structure $J$
compatible with the symplectic structure and consider the space of
$J$-holomorphic spheres $u\colon S^2\to W$
in a given free homotopy class $A\in [S^2,W]$ 
and with the North pole 
on $L_+$ and the South pole on $L_-$. The resulting metric space of
$J$-holomorphic curves is not compact:
the non-compact group of conformal transformations of $S^2$ fixing the 
poles acts properly on it. To make this space compact, we require 
that
\begin{equation}
\labell{eq:normal}
\int_{D_-}u^*\omega=\left<\omega,A\right>/2,
\end{equation}
where $D_-$ is the Southern hemisphere in $S^2$. This condition eliminates
the conformal flow from one pole to another, parametrized by $\R$, and forces
the resulting space $\CM_0$ to be compact. This
space depends on $L_\pm$, $A$, and $J$.

For a generic $J$, the space $\CM_0$ is a smooth manifold
and its free $S^1$-equivariant cobordism class $[\CM_0]$ is independent 
of $J$. (Here the action of $S^1$ on $\CM_0$ arises from the $S^1$-action
on $S^2$ by rotations about the vertical axis.)
Note that in reality
the manifold $W$ and the class $A$ must satisfy some additional conditions,
e.g., $A$ must be minimal, i.e., $\left<\omega,A\right>=m(W,\omega)$;
the definition of $m(W,\omega)$ is recalled below.

Let now $H$ be a smooth non-negative function on $W$ such that $H\equiv 0$
near $L_-$ and $H\equiv \max H$ near $L_+$. Consider the perturbed 
Cauchy--Riemann equation for the maps $u\colon S^2\to W$:
\begin{equation}
\labell{eq:pert-hol}
\bar{\partial}_J u+\lambda \nabla H=0,
\end{equation}
where $\lambda\geq 0$ is real parameter. Note that since
$\nabla H=0$ in some neighborhoods of $L_\pm$, 
this equation simply means that $u$ is holomorphic near the poles.
To make sense of this equation away from the poles, we view $S^2$
as the cylinder $S^1\times \R$ with the two poles attached. Then
the coordinates on the cylinder are used to define 
$\bar{\partial}_J u$ as a vector field along $u$. 

Let $\CM_\lambda$ be the space of solutions of \eqref{eq:pert-hol}
which are in the class $A$, send the North and South poles to $L_\pm$,
and satisfy \eqref{eq:normal}. When $\lambda=0$, we obtain the
space $\CM_0$ introduced above.

In general, the solutions of \eqref{eq:pert-hol} can be thought of
as gradient trajectories for $A_{\lambda H}$ connecting trivial periodic orbits
on $L_+$ with those on $L_-$. Then, a calculation shows that the
difference of actions on these periodic orbits, which is $\lambda\max H$,
is bounded from above by the symplectic area of $u$, i.e.,
$$
\lambda \max H\leq \left<\omega,A\right>. 
$$
As a consequence, $\CM_\lambda=\emptyset$ when $\lambda$ is large.

Let us examine now the disjoint union of the spaces $\CM_\lambda$. 
Under suitable genericity hypotheses, this is a smooth manifold. If this
manifold is compact, it gives a cobordism from $\CM_0$ to the empty set, i.e.,
$[\CM_0]=0$. Just as in Floer's theory, compactness can fail only when
a family of solutions $u_\lambda\in\CM_\lambda$ converges as
$\lambda\to\lambda_0\leq\left<\omega,A\right>/\max H$ to a broken
solution which ``hangs up'' on a contractible one-periodic orbit of 
$\lambda_0 H$, i.e., a contractible $1/\lambda_0$-periodic orbit of $H$.
We conclude that $H$ must have contractible one-periodic orbits,
provided that $\max H> \left<\omega,A\right>$ and $[\CM_0]\neq 0$.
(Note that at this point one still has to show that the orbits found 
are non-trivial.) 

This approach leads to a few versions of the almost existence theorem,
all relying on the same basic requirement that the space $\CM_0$ is 
not cobordant to zero. In other words, the space of holomorphic curves 
``connecting'' $L_+$ and $L_-$ should be sufficiently large for 
the method to apply. In particular, it may be helpful to start
with larger submanifolds $L_\pm$ or with a function $H$ such that
$H\equiv 0$ and $H\equiv \max H$ on large subsets. (The trade-off is that
this may result in the almost existence theorem for a restricted class
of functions $H$.) Let us illustrate these considerations by some examples.

When $L_\pm$ are points and the method is applicable to $W$ and $A$, 
we conclude that $\CHZ(W)\leq \left<\omega,A\right>$ and hence prove
the almost existence theorem in $W$. There are, however, rather few manifolds 
$W$ for which $[\CM_0]\neq 0$ when $L_\pm$ are points. One example is $\CP^n$ 
with the standard symplectic
form normalized as the reduction of the unit sphere in $\C^{n+1}$. In this
case the standard $J$ is already generic, $\CM_0=S^1$ with $S^1$ acting by
translations, and we see that $\CHZ(\CP^n)=\pi$, \cite{HV2}. It is also
possible that this reasoning can be utilized to show that other coadjoint 
orbits of compact semi-simple Lie groups have finite Hofer--Zehnder capacity.

Let now $W=P\times S^2$, where $(P,\beta)$ is a compact symplectic manifold,
as in \cite{HV2}. Then taking the class of the fiber $S^2$  as $A$ and 
applying, with some modifications, this method to $L_-=P\times\{0\}$ and  
$L_+=P\times\{\infty\}$, one can
show that the relative capacity 
$\CHZ(P\times D^2_R, P\times\{0\})=\pi R^2$ as long as
$\area (D^2_R)=\pi R^2<m(P,\beta,J)$, where
$$
m(P,\beta,J):=\inf\{\int_{S^2}u^*\beta>0\mid
\text{ $u$ is a non-constant $J$-holomorphic sphere}\}.
$$
In particular, $m(P,\beta,J)=\infty$ when $\pi_2(P)=0$, \cite{FHV}.
Note also that $m(P,\beta,J)$ is always positive by Gromov's
compactness theorem. Furthermore, assume that
$$
m(P,\beta):=\inf\{\int_{S^2}u^*\beta>0\mid
u\colon S^2\to P\} \geq 0.
$$
In contrast with $m(P,\beta,J)\geq m(P,\beta)$, this constant can be zero. 
(For instance,
this is the case for $P=S^2\times S^2$ where the area of the first component
is 1 and the second component has an irrational area.)  By taking a 
point in $P\times\{0\}$ as $L_-$ and  $L_+$ as before, one can show that 
$\CHZ(P\times S^2,L_+)=\area(S^2)$ as long as $\area(S^2)<m(P,\beta)$,
\cite{HV2}. As a consequence, $\CHZ(P\times D^2_R)=\pi R^2$ if
$\pi R^2 <m(P,\beta)$. (In fact, there are strong indications that
$\CHZ(P\times D^2_R)=\pi R^2$ for any
geometrically bounded symplectic manifold $P$ and any $R>0$; see
\cite{MDS,FS3}.) Symplectic manifolds
$P$ for which $m(P,\beta)>0$ are called rational.

The method also extends to the setting where $P$ is not closed but
is geometrically bounded, \cite{Lu}. (These results have further applications
to the existence problem for periodic orbits of a charge in a magnetic field,
see \cite{Ma,mac2}.) For some other incarnations and applications
of the holomorphic curve 
method, we refer the reader to 
\cite{Lalond,Ke:new,LT,Lu:uniruled,Lu:GW,Lu:toric}. 
In particular, for uniruled symplectic manifolds (such as $\T^2\times S^2$)
and symplectic toric manifolds the Weinstein conjecture was established in
\cite{Lu:uniruled,Lu:toric}.
In the context of the generalized Weinstein--Moser conjecture, it was
proved in \cite{Ke:new} that a small neighborhood of a rational symplectic 
submanifold $M\subset W$ has finite Hofer--Zehnder capacity $\CW$.

Finally, one may replace holomorphic spheres by holomorphic curves
of higher genus to  obtain a sufficiently large space $\CM_0$
as in \cite{LT,Lu:GW}. Extra care in interpreting 
\eqref{eq:pert-hol} is then needed. For instance,
imposing some additional conditions on $H$ may be necessary, e.g., requiring
$H$ to be locally constant outside of a shell separating $L_-$ from 
$L_+$. This, depending
on the details of the approach, leads to a version of either the
nearby existence or almost existence theorem, cf. 
\cite{LT,Lu:uniruled,Lu:GW,Lu:toric}.

The holomorphic curve approach  in the form outlined above 
does not apply to 
compact manifolds that generically have too few holomorphic curves. Among 
such manifolds is, for example,  $\T^4$ with the standard symplectic
structure. It is not known
whether or not the Hofer--Zehnder capacity of $\T^4$ is finite.

\subsection{Hofer's geometry and almost existence}
\labell{sec:Hofer-appl}
Applications of Hofer's geometry to nearby and almost existence theorems
are based on a principle relating minimizing properties of geodesics in
Hofer's metric and fast periodic orbits. Namely, consider the Hamiltonian
flow $\varphi_H^t$ which is a one-parameter subgroup in 
$\CalD$ and hence can be viewed as a geodesic in $\CalD$. Then, 
conjecturally, $\varphi_H^t$ is length minimizing for $t\in [0,\,1)$, 
provided that $H$ has no fast periodic orbits; see \cite{MDS} and references
therein. Various particular cases of this conjecture have been proved. 
We refer the reader to \cite{polt:book} for a detailed discussion 
and additional references and to \cite{MDS}, following \cite{LMD2}, for 
more recent results relevant to our discussion.

The application of Hofer's geometry to the circle of questions considered in
this paper was pioneered in \cite{pol2}.
Below we will just briefly indicate the logic of the argument following
closely
the most recent work \cite{FS3} and in fact suppressing the connection with
minimizing properties of geodesics in Hofer's metric. 

Let $W$ be a geometrically bounded symplectic manifold, which we do not
require to symplectically aspherical or convex, and let
$H$ be an autonomous Hamiltonian supported in $U\subset W$. The
cornerstone of the method is the inequality
\begin{equation}
\labell{eq:cap-displace}
\CW(U)\leq 4e(U),
\end{equation}
which implies almost existence in $U$, provided that $e(U)<\infty$.
(Note that stronger inequalities \eqref{eq:cap-displace1} and 
\eqref{eq:cap-displace1a} hold when $W$ is symplectically aspherical and 
convex. Hence,
the emphasis here is on extending the class of manifolds for which
a version of the capacity--displacement energy inequality is proved.) 

The proof of \eqref{eq:cap-displace} is based on the following two results:
\begin{enumerate}

\item[(i)] Assume that $\supp H \subset U\subset W$ and
$\parallel H\parallel_{\H}>4e(U)$. Then 
$\rho(\varphi_H)<\parallel H\parallel_{\H}$.

\item[(ii)] $H$ has contractible fast periodic orbits, provided that
$\rho(\varphi_H)<\parallel H\parallel_{\H}$.

\end{enumerate}
The first assertion (i) is proved by the curve shortening method, ubiquitous
in Hofer's geometry and going back to \cite{Sic}; see the references 
above for other incarnations of this method. As stated, the second assertion 
(ii) is obtained in \cite{FS3} as a consequence of the results of \cite{MDS}.
(This is the point where length minimizing geodesics enter the picture.)

The inequality \eqref{eq:cap-displace} is sufficient to establish the
almost existence theorem for only a very limited class of domains $U$. 
However, \eqref{eq:cap-displace} can be combined with Macarini's 
stabilization trick, \cite{Ma}. This leads to an upper bound similar
to \eqref{eq:cap-displace}, but with $e(U)$ replaced by the displacement
energy of $U\times S^1$ in $W\times T^*S^1$, \cite{FS3}, which has a 
somewhat broader range of applications. In particular, by
modifying \eqref{eq:cap-displace} in this way, one can prove
the almost existence theorem for low energy periodic orbits of
a charge in a non-vanishing magnetic field and a conservative force field, 
\cite{FS3}. 

We conclude this section by pointing out that it may be possible to prove
\eqref{eq:cap-displace} directly, without invoking minimizing
properties of geodesics in Hofer's geometry, by adapting the argument from
\cite{MDS}.

\section{The Hofer--Zehnder capacity function}
\labell{sec:cap-fun}
What we see as a central open problem concerning the almost existence theorem
is the question whether or not this theorem is sharp. Consider, for
example, a smooth proper  function $H\colon\R^{2n}\to\R$ which we assume to be
bounded from below. By the 
almost existence theorem,  almost all regular levels in $[\min H,\infty)$ carry
periodic orbits. The counterexamples to the Hamiltonian Seifert conjecture
show that $H$ may have a discrete collection of aperiodic levels, i.e.,
regular levels without
periodic orbits, \cite{gi:seifert95,gi:seifert97,gi:bayarea,gg1,gg2,He:fax,He,ke2}.
Moreover, such levels can accumulate to a degenerate critical level of $H$,
\cite{gg3}. However, it is still unknown if aperiodic levels
can accumulate to a regular level either with or without periodic orbits.
In particular, it is not known whether the set of aperiodic energy values can
be dense or be a Cantor set.

\subsection{The definition of the capacity function}
Throughout this section, we will focus on functions on $\R^{2n}$ although
most of our discussion carries over to functions on any manifold of bounded 
capacity. To concentrate on the essential part of the problem, let us assume
that the function has only one critical point, the minimum. The key to
the proof of the almost existence theorem is the Hofer--Zehnder
capacity function associated with $H$:
$$
\cf_H(h)=\CHZ(\{H<h\}) < \infty,\text{ where $h>\min H$.}
$$
This is a monotone increasing function on $(\min H,\infty)$
and the levels $H=h$, where
$\cf_H$ is Lipschitz, carry periodic orbits, \cite{HZ:book}. (This is a 
simple consequence of the Arzela--Ascoli theorem and shortly
we will recall the argument.) Since $\cf_H$ is monotone increasing, it
is differentiable almost everywhere and the almost existence theorem follows.
(Note that there is no reason to expect every non-Lipschitz value
to be aperiodic. In particular, aperiodic points are not entirely visible
from the properties of $\cf_H$.)

Any zero measure set is the set of non-Lipschitz points of some monotone
increasing function, \cite[p. 214]{Na}. Hence, one possible approach to
the problem (but not the only one) is to investigate additional
properties of the capacity function which distinguish it from 
an arbitrary monotone
increasing function. We will soon see that some of these readily arise
from the Arzela--Ascoli theorem.

The following elementary observation illustrates our point, 
\cite{gi:barcelona}. Denote by $Z$ the collection of aperiodic values of $H$.
Furthermore, let $Z_T$ be the collection of levels where all periodic orbits
have period greater than $T$. It is clear that
$$
Z=\bigcap_{T\in\Z_+} Z_T
$$
and that, by the Arzela--Ascoli theorem, $Z_T$ is open. Hence, $Z$ is a 
$G_\delta$ set. Furthermore, if $Z$ is dense, every set $Z_T$ is also
dense. From this we infer that $Z$ must be a residual set when $Z$ is dense.
In particular, the set $Z$ of aperiodic values cannot be a countable dense
set.

\subsection{The area--period relation} The role of the capacity function
in the proof of the almost existence theorem can be best understood in
terms of the classical area--period relation; see, e.g.,
\cite[p. 282]{Ar}. Let us recall this result.

Let $H$ be a function on $\R^2$. Denote by $\area (h)$ the
area bounded by a regular level $H=h$ and by $T(h)$ the period of
the periodic orbit on this level. (Here we are assuming that the levels are
connected.)

\begin{Proposition}[Area--period relation]
\labell{prop:apr}
$d \area (h)/ dh= T(h)$.
\end{Proposition}

Note that the same is true for any proper function on a symplectic
surface, provided that the level $H=h$ is regular. When this level 
is comprised of more than one connected component, the right hand side
is the sum of their periods.

To prove Proposition \ref{prop:apr}, denote by $\omega$
the area form (i.e., the symplectic form) on $\R^2$. 
Dividing $\omega$ by $dH$ near the level, we can write 
$\omega=dH\wedge\alpha$, where $\alpha$ is a one-form such that 
$\alpha(X_H)=1$. Let $\Pi_\eps$ be the annulus $h\leq H\leq h+\eps$. 
The following calculation concludes the proof:
$$
\frac{d\area(h)}{dh}
=\lim_{\eps\to 0}\frac{1}{\eps}\int_{\Pi_\eps} \omega
=\lim_{\eps\to 0}\frac{1}{\eps}\int_{\Pi_\eps} dH\wedge \alpha
=\int_{\{H=h\}}\alpha = T(h).
$$

This result generalizes to higher dimensions when $H$ is convex.
Namely, in this case, on any level $H=h$ there exist periodic orbits 
$\gamma_l$ and $\gamma_r$  such that the periods of these orbits are
equal to the left and, respectively, right derivatives of $\cf_H$ at $h$
and the symplectic areas bounded by $\gamma_l$ and $\gamma_r$ are
equal to $\cf_H(h)$, \cite{Ne}. When the convexity assumption is dropped,
the derivative of the capacity function gives only an upper bound on the
period:

\begin{Proposition}
\labell{prop:apr2}
Assume that the lower derivative 
$$
\LD \cf_H(h) =\liminf_{\delta\to 0}
\frac{\cf_H(h+\delta)-\cf_H(h)}{\delta}
$$ 
is finite. Then the level $H=h$ carries a periodic orbit with period 
$T\leq \LD \cf_H(h)$ and, as a consequence, $\LD \cf_H(h)>0$.
\end{Proposition}

The proposition shows, in particular, that the Hofer--Zehnder capacity
is strictly increasing. More specifically, we have

\begin{Corollary}~
\begin{enumerate}
\item[(i)] The capacity function $\cf_H$ is strictly increasing.

\item[(ii)] Let $U$ and $V$ be open bounded subsets of $\R^{2n}$ such that
$\overline{U}\subset V$. Then $\CHZ(U)<\CHZ(V)$.

\end{enumerate}
\end{Corollary}

It is easy to see that the assumption that the sets are bounded is
essential and the requirement $\overline{U}\subset V$ cannot be replaced
by $U\subsetneq V$.

\begin{proof}[Proof of Proposition \ref{prop:apr2}]
The proof of the proposition is a modification of an argument from 
\cite{HZ:book}. Let $h_i\to h$ be a sequence such that
$$
A_i:=\frac{\cf_H(h_i)-\cf_H(h)}{h_i-h}\to A,
$$
where $0\leq A<\infty$. By passing if necessary to a subsequence, we may
assume that $h_i$ converges to $h$ either from the right or from the left.
In what follows we will assume that $h_i>h$. The other case can be handled
in a similar fashion.
Observe that either $A_i>0$ for all (sufficiently large) indexes $i$ or $\cf_H$
is constant on some interval $[h,\,h+\delta]$, with $\delta>0$, 
since $\cf_H$ is monotone.

Assume first that $A_i>0$ and pick a sequence $b_i>1$ converging to one and 
a sequence 
$$
0<\eps_i< (b_i-1)(h_i-h)A_i.
$$ 
By the definition of the Hofer--Zehnder
capacity, on the domain $\{H<h\}$ there exists an admissible function $K_i$ 
without fast non-trivial periodic orbits and such that 
$\max K_i=\cf_H(h)-\eps_i$. Let $f_i$
be a monotone decreasing function on the interval $[h,\,h_i]$ identically 
equal to $b_i(h_i-h)A_i$ near $h$ and zero near $h_i$ and such that
$|f_i'|\leq a_ib_iA_i$ for some sequence $a_i\to 1+$.
Set $F_i=f_i\circ H$ on the shell $h\leq H\leq h_i$ and smoothly extend this
function to $\R^{2n}$ by requiring it to be constant inside and outside
of the shell. Then
$$
\max (K_i+F_i)=\cf_H(h_i)+[(b_i-1)(h_i-h)A_i-\eps_i]>\cf_H(h_i).
$$
Hence, $K_i+F_i$ must have a non-trivial fast periodic orbit which
can only be a periodic orbit of $F_i$. Therefore, $H$ has a periodic
orbit with period less than or equal to $a_ib_iA_i$ in the shell
$h\leq H\leq h_i$. By passing to the limit and applying the Arzela--Ascoli
theorem, we conclude that $H$ has a periodic orbit of period less than
or equal to $A=\lim a_ib_i A_i$ on the level $H=h$.

To finish the proof it suffices to show that $\cf_H$ cannot be constant 
on any interval $[h,\,h+\delta]$. (This will also be a direct proof of the
corollary.) Assume the contrary. Then $\cf_H(h_i)=\cf_H(h)$
for all large enough $i$. In this case we again choose a sequence of
decreasing functions $f_i$ on $[h,\,h_i]$ such 
that $f_i$ is zero near $h_i$ and positive constant near $h$. (Thus
$f_i$ is equal to $\max f_i>0$ near $h$.)
Furthermore, we can make the slope of $f_i$ so small that the
function $F_i$ has no fast non-trivial periodic orbits.
Let $0<\eps_i<\max f_i$ and let $K_i$ be as above. It is clear that
the function $K_i+F_i$ has no fast non-trivial periodic orbits. On the
other hand, we have
$$
\max (K_i+F_i)=\cf_H(h)-\eps_i+\max f_i> \cf_H(h)=\cf_H(h_i),
$$
and hence $K_i+F_i$ must have fast non-trivial periodic orbits. This 
contradiction completes the proof.
\end{proof}

A question closely related to the area--period relation is that of
representability of a capacity. Let $U$ be the domain bounded
by a compact smooth hypersurface $\Sigma$ in $\R^{2n}$. We say that
a capacity $\cf$ is represented on $\Sigma$ if there is a
closed characteristic on $\Sigma$ that bounds a disc of symplectic area 
$\cf(U)$.
For example, as is easy to see from the definition, $\wsh$ is represented
on smooth hypersurfaces of contact type. Since there are hypersurfaces
without closed characteristics, no capacity is represented on every
hypersurface. The Hofer--Zehnder capacity is represented on convex
hypersurfaces, \cite{HZ:cap}. In fact, when $\Sigma$ is convex,
$\CHZ(U)$ is the minimal symplectic area of a closed characteristic
on $\Sigma$, \cite{HZ:cap}. However, it is not known if $\CHZ$ is
represented on every contact type, or even restricted contact type,
hypersurface. Arguing as in the proof of Proposition \ref{prop:apr2}
and utilizing Lemma \ref{lemma:sigma}, one can prove the following

\begin{Proposition}
\labell{prop:repr}
Let $U$ be the domain bounded by a smooth hypersurface $\Sigma$ of
restricted contact type. Then, $\CHZ$ is sub-represented on $\Sigma$,
i.e., there exists a closed characteristic
on $\Sigma$ with area less than or equal to $\CHZ(U)$. 
\end{Proposition}

Recall in this connection that even when $U$ is star-shaped there can be
a closed characteristic on $\Sigma$ with area strictly less than $\CHZ(U)$;
the ``Bordeaux bottle'' (see \cite[p. 99]{HZ:book}) is an example.
We refer the reader to, e.g., \cite{FS1,Her2} for other representability
results.

\subsection{Period growth and the Hofer--Zehnder capacity function}
As is immediately clear from the definition, the Hofer--Zehnder capacity 
function is 
necessarily continuous from the left. Beyond this trivial observation,
little is known about continuity or differentiability properties of
the  capacity function. For example, it is not known if
this function can be discontinuous at regular values. One, admittedly 
very naive, approach
to this question is based on the estimates of the period growth.

To make this more precise, fix a regular level of $H$, say $H=0$. Let
$\tau(h)$ be the infimum of all periods of periodic orbits in the
open shell $0<H<h$ or $h<H<0$, depending on whether $h$ is positive or
negative. Then the growth of the function $\tau(h)$ is related to
continuity and smoothness of $\cf_H$ at zero as the next proposition
shows.

\begin{Proposition}
\labell{prop:growth}~
\begin{enumerate}
\item[(i)] $\tau(h)=O(1/h)$ as $h\to 0+$ and $\tau(h)=o(1/|h|)$ as $h\to 0-$.

\item[(ii)] Assume that $\tau(h_k)\geq C/h_k$ for some sequence $h_k\to 0+$.
Then $\cf_H(0+)-\cf_H(0)\geq C$.

\item[(iii)] Assume that $\tau(h_k)\geq C/|h_k|^\alpha$ for some sequence 
$h_k\to 0$ and $0\leq\alpha\leq 1$. Then 
$|\cf_H(h)-\cf_H(0)|\geq C|h|^{1-\alpha}$.
\end{enumerate}
\end{Proposition}

The proofs of these facts follow the same line as the proof of the
almost existence theorem in \cite{HZ:book}, cf. the proof of Proposition
\ref{prop:apr2}. Omitting a detailed argument here,
we only mention that (i) is a consequence of the fact that
the capacity is bounded. In fact, $\CHZ(\{H<0\})=\infty$ 
if (i) fails for the left limit and $\CHZ(V)=\infty$ for any open set
$V\supset \{H\leq 0\}$ if (i) fails for the right limit.

\begin{Remark}
The analogues of (ii) and (iii) hold when the difference of 
capacities is replaced by the relative capacity.
\end{Remark}

Proposition \ref{prop:growth} is difficult to apply to 
examine discontinuity or non-smoothness of the capacity function.
Indeed, let, for example, $H$ be a Hamiltonian on $\R^{2n}$ constructed in 
\cite{gi:seifert95,gi:seifert97,gg1,gg2,He:fax,He,ke2}
such that $H=0$ carries no periodic orbits. Then $\cf_H$ is not
Lipschitz at zero as is clear already from the results of \cite{HZ:book}.
To utilize Proposition \ref{prop:growth}, we would need to bound from
below minimal
periods on the nearby levels $H=h$. However, the constructions of $H$ afford
little insight into the dynamics on these levels and obtaining such lower
bounds appears to be an extremely difficult problem. Preliminary estimates 
indicate
that $\cf_H$ can probably be H\"older with any $\alpha\in (0,\,1)$ at $H=0$
or even fail to be H\"older. However, there is no convincing evidence that
$\cf_H$ can be discontinuous. (Note in this connection that the property 
of $\cf_H$ to be H\"older with
a specific $\alpha$ or to be continuous, just as to be Lipschitz 
(see \cite{HZ:book}), is 
determined by the level $H=0$ and is independent of the choice of $H$.)

Consider, however, the following example. 

\begin{Example}[The horocycle flow]
\labell{ex:horocycle}
Let $M$ be a closed surface equipped
with a metric of constant negative curvature $-1$ and let $\Omega$ be
the area form on $M$. Consider the twisted symplectic structure 
$\omega=\omega_0+\pi^*\Omega$ on $W=T^*M$, where $\pi$ is the
natural projection $T^*M\to M$ and $\omega_0$ is the standard symplectic
structure. Set $H=\parallel p\parallel^2-1$. The Hamiltonian flow on the 
level $H=0$ is the horocycle flow and hence has no periodic orbits; see, e.g., 
\cite{gi:survey,gi:bayarea} for details. Assume that the Hofer--Zehnder
capacity of the sets $\{H<h\}$ for $h$ near zero is finite. (This is not
known, although is likely to be true, even for $h<0$ close to zero. For
small neighborhoods of the zero section (i.e., $h\approx -1$),
finiteness has been recently proved in \cite{FS1}; however, for
$h>0$, finiteness of $\CHZ$ appears to be beyond reach of the methods 
considered in this paper.) Then, we claim that
\begin{equation}
\labell{eq:horocycle}
\cf_H(h)-\cf_H(0)\geq l_{\min}\sqrt{h} \text{ for $h>0$},
\end{equation}
where $l_{\min}$ is the minimal length of a closed geodesic on $M$.
Therefore, $\cf_H$ is not smoother at $h=0$ than $1/2$-H\"older on the right.

To prove this, let us first note that
\begin{equation}
\labell{eq:hard}
\tau(h)=l_{\min}/\sqrt{h}\text{ for $h>0$}, 
\end{equation}
 and hence, by (iii),
$\cf_H(h)-\cf_H(0)\geq l_{\min}\sqrt{h}$. Establishing \eqref{eq:hard}
directly appears to be rather difficult. Instead, let us
observe that there exists a symplectomorphism 
$(T^*M\ssminus \{H\leq 0\},\omega)\to (T^*M\ssminus M, \omega_0)$ which
sends $\parallel p\parallel^2$ to $\parallel p\parallel^2-1$. (Such a 
symplectomorphism can be constructed, for instance,  by judiciously applying 
Moser's method.) This allows one to translate 
the period growth on $(T^*M,\omega_0)$ as $\parallel p\parallel\to 0$, 
which is $l_{\min}/\parallel p\parallel$
for the Hamiltonian $\parallel p\parallel^2/2$, to the period growth on 
$(T^*M,\omega)$ for $H>0$.

Note also that by applying the results of \cite{BPS} and using
the above symplectomorphism, one can prove that 
$\CHZ(\{ H<h\},\{H\leq 0\})\leq l_{\min}\sqrt{h}$. This,
together with the analogue of \eqref{eq:horocycle} for the relative
capacity, proves that $\CHZ(\{ H<h\},\{H\leq 0\})= l_{\min}\sqrt{h}$.

One can also estimate the period growth on the left. All
orbits on the level $H=h<0$ are closed and project to geodesic circles
on $M$ with geodesic curvature $k_g>1$. A straightforward calculation 
shows that $\tau(h)=2\pi/\sqrt{|h|}$ for $h<0$ and hence, by (iii),
$\cf_H(0)-\cf_H(h)\geq 2\pi\sqrt{|h|}$. Thus, on the left, $\cf_H$ is again
not smoother than $1/2$-H\"older.
\end{Example}

\end{document}